\documentstyle[12pt,amscd]{amsart}
\newtheorem{Theorem}{Theorem}[section] \newtheorem{Lemma}[Theorem]{Lemma}
\newtheorem{Corollary}[Theorem]{Corollary}
\newtheorem{Proposition}[Theorem]{Proposition}

 \def\depth{\operatorname{depth}}
\def\reg{\operatorname{reg}} \def\ini{\operatorname{in}}
\def\Gin{\operatorname{Gin}}
\def\Ker{\operatorname{Ker}} \def\To{\longrightarrow}
\def\Ext{\operatorname{Ext}} 
 \def\bideg{\operatorname{bideg}}

\def\To{\longrightarrow}
\def\sk{\smallskip}
\def\mm{{\frak m}}
\def\nn{{\frak n}}
\let\union=\cup

\begin{document}
\title{Gr\"obner bases and regularity of Rees algebras}
\author{J\"urgen Herzog, Dorin Popescu and Ng\^o Vi\^et Trung}
\address{Fachbereich Mathematik,  Universit\"at-GHS Essen, 45117 Essen, Germany}
\email{juergen.herzog@@uni-essen.de}
\address{Institute of Mathematics, University of Bucharest, P.O. Box 1-764,
Bucharest 70700, Romania} \email{dorin@@stoilow.imar.ro}
\address{Institute of Mathematics,  Box 631, B\`o H\^o, 10000 Hanoi, Vietnam}
\email{nvtrung@@hn.vnn.vn}
\thanks{The third author is partially supported by the National Basic
Research Program of Vietnam} \keywords{} \subjclass{} \maketitle

\section*{Introduction} \smallskip

Let $B=k[x_1,\ldots,x_n]$ be a polynomial ring over a field $k$
and $A=B/J$ a quotient ring of $B$ by a homogeneous ideal $J$.
Let $\mm$ denote the maximal graded ideal of $A$. Then the Rees
algebra $R=A[\mm t]$ may be considered a standard graded
$k$-algebra and has a presentation $B[y_1,\ldots,y_n]/I_J$. For
instance, if $A = k[x_1,\ldots,x_n]$, then
$$R \cong k[x_1,\ldots,x_n,y_1,\ldots,y_n]/(x_iy_j-x_jy_i|\ i, j =
1,\ldots,n).$$ In this paper we want to compare the ideals $J$ and
$I_J$ as well as their homological properties.\sk

The generators of $I_J$ can be easily described as follows. For
any homogeneous form $f=\sum_{1\leq i_1\leq\cdots\leq i_d\leq
n}a_{i_1\cdots i_d}x_{i_1}\cdots x_{i_d}\in B$ of degree $d$ we
set
\[
f^{(k)} := \sum_{1\leq i_1\leq\cdots\leq i_d\leq n}a_{i_1\cdots
i_d}x_{i_1}\cdots x_{i_{d-k}}y_{i_{d-k+1}}\cdots y_{i_d}
\]
for $k=0,\ldots,d$. For any subset $L \subset B$ of homogeneous
polynomials in $S$ we set
$$L' := \{f^{(k)}|\ f\in L, k=0,\ldots, \deg f\},$$
and let $$H :=\{x_iy_j-x_jy_i|\ 1\leq i<j\leq n\}.$$ If $L$ is a
minimal system of generators of $I$, then $L'\union H$ is a
minimal system of generators of $I_J$ (Proposition \ref{kernel}).
We will show that if $L$ is Gr\"obner basis of $J$ for the
reverse lexicographic order induced by $x_1>\cdots>
x_n>y_1>\cdots >y_n$, then $L'\union H$ is Gr\"obner basis of
$I_J$ (Theorem \ref{groebner}). As a consequence, if $J$ has a
quadratic Gr\"obner basis, then so does $I_J$. \sk

The main concern of this paper is however the regularity which is
a measure for the complexity of the resolution of a standard
graded algebra (see [EiG], [BM]). Recall that the
Castelnuovo-Mumford regularity of $A$ is defined by
$$\reg(A):= \max\{b_i-i|\ i >
0\},$$ where $b_i$ denotes the largest degree of a generator of
the $i$th syzygy module of $A$. The regularity and related
invariants of a graded $k$-algebra (for example, the extremal
Betti numbers introduced in [BCP]) can be expressed in terms of
the cohomological invariants $a_i=\max\{a|\ H^i(A)_a\neq 0\}$,
where $H^i(A)$ denotes the $i$th local cohomology of $A$ with
support $\mm$ (see Section 2 for more details). For instance,
$\reg(R) = \max\{a_i +i|\ i \ge 0\}$. In particular, we will
also study the invariant
$$a^*(A) := \max\{a_i|\ i \ge 0\} = \max\{b_i|\ i \ge 0\} - n,$$
which is another kind of regularity for $A$ [Sh], [T2], [T3]. Our
results are based on the observation that the local cohomology of
$R$ can be estimated in terms of the local cohomology of  $A$. \sk

If $R$ is the Rees algebra of an arbitrary homogeneous ideal $I$
of $A$ generated by forms of the same degree, then $R$ is still a
standard $k$-algebra. In this case, we have the following
estimations:
\begin{eqnarray*} a^*(A)-s & \le & a^*(R)\ \le\
\max\{a^*(A),a^*(G)\},\\ \reg(A) & \le & \reg(R)\ \le\
\max\{\reg(A)+1,\reg(G)\}, \end{eqnarray*} where $s$ is the
minimal number of generators of $I$ and $G$ denotes the associated
graded ring of $I$ (Theorem \ref{a-invariant1} and Theorem
\ref{regularity1}). These bounds are sharp. In particular, if $R$
is the Rees algebra of the graded maximal ideal of $A$, then
\begin{eqnarray*} a^*(A) -n & \le & a^*(R)\ \le\ a^*(A),\\
\reg(A) & \le & \reg(R)\ \le\ \reg(A)+1. \end{eqnarray*}  It is
shown in Theorem \ref{a-invariant3} and Theorem \ref{regularity3}
that $a^*(R) = a^*(A)$ if and only if $a^*(A) \neq -1$ and that
$\reg(R)= \reg(A)+1$ if and only if there is an integer $i$ such
that $\reg(A)=a_i+i$ and $a_i\leq -2$. The proofs follow from the
fact that the bigraded components of the local cohomology of $R$
can be expressed completely by the graded components of the local
cohomology of $A$ (see Theorem \ref{vanish3}). In particular, we
can show that $\reg(R)=\reg(A)+1$ if $\reg(A)=b_i-i$ and $b_i\leq
n-2$ for some index $i$ at which $A$ has an extremal Betti number
(Corollary \ref{sufficient}). However, an example shows that this
condition is only sufficient. As applications, we compare the
regularity of the Rees algebra of the ring $B/\ini(I)$, where
$\ini(I)$ denotes the initial ideal of $I$, with that of $R$ and
we estimate this regularity for the generic initial ideal
$\Gin(J)$ with respect to the reverse lexicographic term order.
\sk

We will also compute the projective dimension of $I_J$. In
Proposition \ref{depth} we give a precise formula for the depth
of $R$ in terms of invariants of $A$. In fact,
$$\depth R = \max\{i|\ H_{\frak m}^j(A)_a = 0\ \text{\rm for}\ a \neq
-1, j < i-1,\ \text{\rm and}\ a_{i-1} < 0\}.$$  This formula is
better than Huckaba and Marley's estimation for the depth of the
Rees algebra of an arbitrary ideal in a local ring [HM]. Inspired
by a construction of Goto [G] we give examples showing that for
arbitrary positive numbers $2\leq r<d$ there exists a standard
graded $k$-algebra $A$ of dimension $d$ with $\depth A=r$ and
$\depth R=d+1$. In these examples $R$ is Cohen-Macaulay, since
$\dim R=d+1$. Though the difference between the depth of $A$ and
of $R$ may be large, this is not the case for the Rees ring $R^*$
of a polynomial ring extension $A[z]$ of $A$. Here we have that
$\depth R^*=\depth A+1$ if $a_s\geq 0$ and $\depth A[z]=\depth
A+2$ if $a_s<0$ (Corollary \ref{difference}). \sk

We would like to mention that if $R$ is the Rees algebra of a
homogeneous ideal generated by forms of different degree, $R$ is
not a standard $k$-algebra. Since $R$ is a standard graded
algebra over $A$, one can still define the Castelnuovo-Mumford
regularity and the $a^*$-invariant of $R$ with respect to this
grading. These invariants have been studied recently by several
authors (see e.g. [JK], [Sh], [T1], [T2]). \sk

\noindent{\it Acknowledgement.} This paper was completed during a
visit of Ng\^o Vi\^et Trung to the University of Essen im summer
2000. He would like to express his sincere thanks to the Research
Group "Arithmetic and Geometry" and to J.~Herzog for support and
hospitality.\sk

\section{Gr\"obner basis of Rees algebras} \sk

Let $A$ be a standard graded $k$-algebra with graded maximal
ideal ${\frak m} = (x_1,\ldots,x_n)$. Then  $A=B/J$ where
$B=k[x_1,\ldots, x_n]$ is a polynomial ring, and $J\subset B$ a
graded ideal. The Rees algebra $R=A[{\frak m}t]$ may be
considered as a bigraded module over the bigraded polynomial ring
$S=k[x_1,\ldots,x_n,y_1,\ldots,y_n]$ (where $\deg x_i=(1,0)$ and
$\deg y_i=(1,1)$ for all $i$) via the bigraded epimorphism $\phi:
S\to R$ with $\phi(x_i)=x_i$ and $\phi(y_i)=x_it$ for
$i=1,\ldots,n$. Let $I_J$ denote the kernel of this
epimorphism.\sk

We are interested in the generators and the Gr\"obner basis of
$I_J$. In order to describe $I_J$ we introduce the following
notations. \sk

Let $f=\sum_{1\leq i_1\leq\cdots\leq i_d\leq n}a_{i_1\cdots
i_d}x_{i_1}\cdots x_{i_d}\in B$ be homogeneous of degree $d$. For
$k=0,\ldots,d$ we set
\[
f^{(k)}= \sum_{1\leq i_1\leq\cdots\leq i_d\leq n}a_{i_1\cdots
i_d}x_{i_1}\cdots x_{i_{d-k}}y_{i_{d-k+1}}\cdots y_{i_d}.
\]
Notice that $f^{(k)}$ is bihomogeneous of degree $(d,k)$. For any
subset $L \subset B$ of homogeneous polynomials in $S$ we set
$$L' := \{f^{(k)}|\ f\in L, k=0,\ldots, \deg f\}.$$ We further let
$$H :=\{x_iy_j-x_jy_i|\ 1\leq i<j\leq n\}.$$

With these notations we have\sk

\begin{Proposition}
\label{kernel} Let $L$ be a (minimal) system of generators of
$J$, then $L'\union H$ is a (minimal) system of generators of
$I_J$.
\end{Proposition}

\begin{pf}
Let $P=k[x_1,\ldots,x_n,x_1t,\ldots,x_nt]\subset
k[x_1,\ldots,x_n,t]$, and $\phi_1: S\to P$, $\phi_2: P\to R$ be
the $k$-algebra homomorphisms given by $\phi_1(x_i)=x_i$,
$\phi_1(y_i)=x_it$, and $\phi_2(x_i)=\bar{x}_i$,
$\phi_2(x_it)=\bar{x}_it$ for $i=1,\ldots,n$. We have
$\phi=\phi_2\circ\phi_1$, and since $\phi$ is bigraded, the ideal
$I_J$ is bigraded. We clearly have  $L'\union H \subset I_J$. Let
$f\in I_J$ be bigraded with $\deg f=(a,b)$. Then
$\phi_1(f)=f(x,xt)=f(x,x)t^b$, and so
$0=\phi(f)=f(\bar{x},\bar{x})t^b$, that is,
$f(\bar{x},\bar{x})=0$. Therefore, there exist homogeneous
elements $g_i\in B$ and $f_i\in J$ such that
$f(x,x)=\sum_{i=1}^mg_if_i$. Let $b_i=\min\{\deg f_i,b\}$. Then
\[
\phi_1(f)=f(x,x)t^b=\sum_{i=1}^m(g_it^{b-b_i})(f_it^{b_i})
=\phi_1(\sum_{i=1}^mg_i^{(b-b_i)} f_i^{(b_i)}),
\]
and so $f\in L'\union H$, since $\Ker\phi_1$ is generated by $H$.

Now let $L$ be a minimal system of generators of $J$. We first
show that $\phi_1(L')$ is a minimal system of generators of the
ideal $L_J=\phi_1(I_J)$ in $P$. Indeed, $\phi_1(L')=\{f_it^b|\
f_i\in L,\quad b=0,\ldots, \deg f_i\}$. Suppose this is not a
minimal system of generators of $L_J$. Then there exists an
equation
\[
f_it^b=\sum_{j}\sum_k(f_jt^{b_{jk}})(g_{jk}t^{c_{jk}}),
\]
where $b_{jk}\leq \deg f_j$, $b_{jk}+c_{jk}=b$ and
$f_jt^{b_{jk}}\neq f_it^b$ for all $j$ and $k$, and where all
summands are bihomogeneous of degree $(d, b)$ with $d=\deg f_i$.
Notice that the right hand sum contains no summand of the form
$(f_it^{b_{ik}})(g_{ik}t^{c_{ik}})$. In fact, otherwise we would
have  $\deg g_{ik}t^{c_{ik}} = (0,b-b_{ik})$, and so $b_{ik}=b$
which is impossible. It follows that $f_i=\sum_{j\neq i}(\sum_k
g_{jk})f_j$, a contradiction.\sk

Now suppose that $L'\union H$ is not a minimal system of
generators of $I_J$. If one of the $f_i^{(k)}$ is a linear
combination of the other elements of $L'\union H$, then
$\phi_1(L')$ is not a minimal system of generators of $L_J$, a
contradiction.  Next suppose one of the elements of $H$, say,
$x_1y_2-x_2y_1$  is a linear combination of the other elements of
$L'\union H$. Only the elements of bidegree $(2,1)$ can be
involved in such a linear combination. In other words,
\[
x_1y_2-x_2y_1=\sum \lambda_f f^{(1)}+h\quad\text{with}\quad
\lambda\in k.
\]
Here the sum is taken over all $f\in L$ with $\deg f=2$, and $h$
is a $k$-linear combination of the polynomials $x_iy_j-x_jy_i$
different from $x_1y_2-x_2y_1$.  Since the monomial $x_2y_1$ does
not appear in any polynomial on the right hand side of the
equation, we get a contradiction.
\end{pf}\sk

We will now compute a Gr\"obner basis of $I_J$. For the proof we
will use the following Gr\"obner basis criterion.\sk

\begin{Lemma}
\label{wellknown} Let $Q=k[x_1,\ldots,x_r]$ be the polynomial
ring,  $I \subset Q$ a graded ideal and $L$ a finite subset of
homogeneous elements of $I$. Given a term order $<$, there exists
a unique monomial $k$-basis $C$ of $Q/(\ini(L))$ (which we call a
``standard basis'' with respect to $<$ and $L$). This $k$-basis
$C$ is a system of generators for the $k$-vector space $Q/I$, and
$L$ is a Gr\"obner basis of $I$ with respect to $<$, if and only
if $C$ is a $k$-basis of $Q/I$.
\end{Lemma} \sk

\begin{Theorem}
\label{groebner} Let $<$ be the reverse lexicographic order
induced by $x_1>\cdots> x_n>y_1>\cdots >y_n$. If $L$ is a
Gr\"obner basis of $J$ with respect to the term order $<$
restricted to $B$, then $L'\union H$ is a Gr\"obner basis of
$I_J$ with respect to $<$.
\end{Theorem}

\begin{pf} Let $C$ be a standard basis of $B$ with respect to $<$ and $L$,
and set $$C' :=\{u^{(k)}|\ u\in C,\quad k=0,\ldots,\deg u\}.$$ We
will show that
\begin{enumerate}
\item[(i)] $C'$ is a standard basis with respect to $<$ and
$L'\union H$, and
\item[(ii)] $C'$ is a $k$-basis of $R$.
\end{enumerate}
Let $v$ be a monomial of $T$ which does not belong to the ideal
$(\ini(L')\union\ini(H))$. Since $v\not\in \ (\ini(H))$, it
follows that $v=u^{(k)}$ for some monomial $u\in S$. Suppose that
$u\not\in C$. Then $u\in\ini(L)$, and since $\ini(L')=\ini(L)'$
it follows that $u^{(k)}\in (\ini(L'))$, a contradiction. Thus
$u\in C$, and hence $v\in C'$.

Conversely, if $v\in C'$, then $v=u^{(k)}$ for some $u\in C$.
Monomials of the form $u^{(k)}$ cannot be multiples of monomials
of $\ini(H)$. Suppose $u^{(k)}$ is a multiple of a monomial $w\in
\ini(L)'$. Then $w=v^{(l)}$ for some $v\in\ini(L)$ and some $l$,
and $v^{(l)}$ divides $u^{(k)}$. It follows that $v$ divides $u$,
a contradiction. This proves (i).

Let $C_i=\{u\in C|\ \deg u=i\}$, and similarly $C'_i=\{u\in C'|\
\deg u=i\}$. Since $C$ is a $k$-basis of $A$, it follows that
$|C_i|=\dim_kA_i$, and since $C'_i=\{u^{(k)}|\ u\in C,
k=0,\ldots,i\}$, it follows  that $|C'_i|=(i+1)|C_i|=
(i+1)\dim_kA_i$. It is easy to show that
$\dim_kR_i=(i+1)\dim_kA_i$, $|C'_i|=\dim_kR_i$ for all $i$. This
shows that the elements of $C'$ are $k$-linearly independent, and
proves (ii). Hence the desired conclusion follows from Lemma
\ref{wellknown}.
\end{pf} \sk

\begin{Corollary}
\label{quadratic} If $J$ has a quadratic Gr\"obner basis, then so
does $I_J$.
\end{Corollary}

We would like to remark that if $L$ is a reduced Gr\"obner basis,
then $L'\union H$ need not be reduced.\sk

\noindent{\bf Example.} Let  $A=k[x_1,x_2,x_3]/(x_1x_2-x_3^2)$.
Then $L=\{x_1x_2-x_3^2\}$ is a reduced Gr\"obner basis of $J$,
but $L'\union H$ is not reduced, since
$x_1y_2=\ini(x_1y_2-x_3y_3)$ appears in $x_1y_2-x_2y_1$.\sk

\section{Regularity and local cohomology of graded algebras}\sk

The aim of this section is to prepare some facts on the
relationships between the regularities and local cohomology
modules of a graded module. \sk

Let $B = k[x_1,\ldots,x_n]$ be a polynomial ring over a field
$k$. Let $E$ be a finitely graded module over $B$. Let ${\Bbb F}:
0 \to F_r \to \cdots \to F_1 \to F_0 \to E\to 0$ be a minimal free
resolution of $E$. For all integer $i$ we denote by $b_i(E)$ the
largest degree of the generators of $F_i$, where $b_i(E) :=
-\infty$ for $i < 0$ or $i > r$. The Castelnuovo-Mumford
regularity of $E$ [EiG] is defined by
$$\reg(E) := \max\{b_i(E)-i|\ i \ge 0\}.$$
This notion is refined by D. Bayer, H. Charalambous, and S.
Popescu [BCP] as follows. For any integer $j$ let
$$j\text{-}\reg(E) := \max\{b_i(E)-i|\ i \ge j\}.$$
Similarly, we can define the invariants
\begin{eqnarray*}
b^*(E) & := & \max\{b_i(E)|\ i \ge 0\},\\
b_j^*(E) & := & \max\{b_i(E)|\ i \ge j\}.
\end{eqnarray*}
It is known that these invariants can be also characterized by
means of the graded local cohomology modules of $E$. \sk

Let $A = B/J$ be any graded quotient ring of $B$. Let $E$ now be
a finitely generated module over $A$. Let ${\frak m}$ denote the
maximal graded ideal of $A$. For any integer $i$ we denote by
$H_{\frak m}^i(E)$ the $i$th local cohomology module of $E$. Since
$H_{\frak m}^i(E)$ is a graded artinian $A$-module, $H_{\frak
m}^i(E)_a = 0$ for $a$ large enough. Therefore we can consider
the largest non-vanishing degree $a_i(E) := \max\{a|\ H_{\frak
m}^i(E)_a \neq 0\}$, where $a_i(E) = -\infty$ if $H_{\frak
m}^i(E) = 0$. Note that $H_{\frak m}^i(E) = 0$ for $i < 0$ and $i
> d:= \dim E$ and that $a_d(E)$ is the $a$-invariant of $E$ [GW].
For any integer $j$ we define
\begin{eqnarray*}
a_j^*(E) & := & \max\{a_i(E)|\ i \le j\},\\
\reg_j(E) & := &\max\{a_i(E)+i|\ i \le j\}.
\end{eqnarray*}
In particular, we set
$$a^*(E) :=  \max\{a_i(E)|\ i \ge 0\}.$$ These cohomological invariants do not
depend on the presentation of $A$. See [Sh], [T1], [T2], [T3] for
more information on these invariants.\sk

\begin{Theorem} \label{link1} {\rm [T3, Theorem 3.1]}
For any integer $j$ we have \par \noindent {\rm (i)} $b_j^*(E) =
a_{n-j}^*(E)+n$,\par \noindent {\rm (ii)} $j\text{-}\reg(E) =
\reg_{n-j}(E)$.
\end{Theorem}\sk

Theorem \ref{link1} has the following immediate consequence.\sk

\begin{Corollary} \label{link2}
{\rm (i)} $b^*(E) = \max\{a_i(E)|\ i \ge 0\}+n = a^*(E)+n,$ \par
\noindent {\rm (ii)} $\reg(E) = \max\{a_i(E)+i|\ i \ge 0\}$ {\rm
[EG, Theorem 1.2]}.
\end{Corollary} \sk

>From Theorem \ref{link1} we also obtain the following
relationship between the invariants $b_j(E)$ and $a_{n-j}(E)$.
Following [BCP] we say that $E$ has an {\it extremal Betti number}
at $j$ if  $b_j(E)-j
> b_i(E)-i$ for all $i > j$ or, equivalently, $j\text{-}\reg(E)
> (j+1)\text{-}\reg(E)$. \sk

\begin{Corollary} \label{link3} Assume that $E$ has an extremal Betti
number at $j$. Then $b_j(E) = a_{n-j}(E) + n$. \end{Corollary}

\begin{pf} By the assumption, $j\text{-}\reg(E) = b_j(E)-j$. By Theorem
\ref{link1},
$$\reg_{n-j}(E) = j\text{-}\reg(E) > (j+1)\text{-}\reg(E) =
\reg_{n-j-1}(E).$$ Therefore, $\reg_{n-j}(E) = a_{n-j}(E)+n-j$.
>From this it follows that $b_j(E) = a_{n-j}(E)+n$.
\end{pf}\sk

For later applications we also prepare some facts on the
regularity of polynomial extensions and quotient modules. \sk

\begin{Lemma} \label{extension1} Let $A[z]$ be a polynomial ring over $A$ in
one variable. Let ${\frak n}$ denote the maximal graded ideal of
$A[z]$. Put $E[z] = E \otimes A[z]$ and $E[z^{-1}] = E \otimes
A[z^{-1}]$. For every integer $i \ge 1$ we have
$$H_{\frak n}^i(E[z]) = H_{\frak m}^{i-1}(E)(1)[z^{-1}].$$
\end{Lemma}

\begin{pf} By local duality (see e.g. [BH, Theorem 3.6.19])
we know that
\begin{eqnarray*} H_{\frak m}^i(E) & = &
\Ext_B^{n-i}(E,B(-n))^{\vee},\\ H_{\frak n}^i(E[z]) & = &
\Ext_{B[z]}^{n+1-i}(E[z],B[z](-n-1))^{\vee}, \end{eqnarray*} where
$^{\vee}$ denotes the Matlis duality. Since $B \To B[z]$ is a
flat extension, we have $$\Ext_{B[z]}^{n+1-i}(E[z],B[z]) =
\Ext_B^{n+1-i}(E,B)[z].$$ From this it follows that
\begin{eqnarray*} H_{\frak n}^i(E[z]) & = &
\Ext_B^{n+1-i}(E,B(-n-1))^{\vee}[z^{-1}]\\ & = &
\Ext_B^{n+1-i}(E,B(-n))^{\vee}(1)[z^{-1}]\ =\ H_{\frak
m}^{i-1}(E)(1)[z^{-1}].\end{eqnarray*} \end{pf}\sk

\begin{Proposition} \label{extension2} With the above notation we
have\par \noindent {\rm  (i) } $a_i(E[z]) = a_{i-1}(E) - 1$ for
all $i \ge 0$, \par \noindent {\rm  (ii) } $a^*(E[z]) = a(E)-1$,
\par \noindent {\rm  (iii) } $\reg(E[z]) = \reg(E)$.
\end{Proposition} \sk

\begin{pf} By Lemma \ref{extension1} we have
$$H_{\frak n}^i(E[z])_n = \bigoplus_{a \ge n-1}H_{\frak m}^{i-1}(E)_a.$$
Hence (i) is immediate. The formulas (ii) and (iii) are
consequences of (i). \end{pf} \sk

\begin{Proposition} \label{quotient} Assume that $\depth E > 0$ and
$f \in A $ is a regular form of degree $c$ for $E$. Then \par
\noindent {\rm (i) } $a^*(E/fE) = a^*(E)+c$,\par \noindent {\rm
(ii) } $\reg(E/fE) = \reg(E) + c-1$. \end{Proposition}

\begin{pf} From the exact sequence $0 \To E(-c) \stackrel{f}{\To}
E \To E/fE \To 0$ we obtain the following exact sequence of local
cohomology modules:
$$ H_\mm^i(E)_a \To H_\mm^i(E/fE)_a \To H_\mm^{i+1}(E)_{a-c} \To
H_\mm^{i+1}(E)_a\ .$$ From this it immediately follows that
$a_i(E/fE) \le \max\{a_i(E),a_{i+1}(E)+c\}.$  For $a \ge
a_i(E/fE)$, the map $H_\mm^{i+1}(E)_{a-c} \To H_\mm^{i+1}(E)_a$ is
injective. Since $H_\mm^{i+1}(E)_a = 0$ for all large $a$, this
injective map yields $H_\mm^{i+1}(E)_{a-c} = 0$. Therefore, we get
$a_{i+1}(E)+c \le a_i(E).$  Taking the maxima over $i$ of the
inequalities
$$a_{i+1}(E)+c \le a_i(E/fE) \le \max\{a_i(E),a_{i+1}(E)+c\}$$
we will get (i). For (ii) we only need to take the maxima over $i$
of the inequalities
$$a_{i+1}(E)+ i + c \le a_i(E/fE)+i \le \max\{a_i(E)+i,a_{i+1}(E)+i +
c\}.$$ \end{pf} \sk

\section{Rees algebras of ideals generated by forms of the same degree}
\sk

Let $A$ be a standard graded algebra over a field. Let $I =
(f_1,\ldots,f_s)$ be a homogeneous ideal in $A$ such that
$f_1,\ldots,f_m$ have the same degree. Then the Rees algebra $R =
A[It]$ can be considered as a standard $\Bbb N$-graded algebra
over $k$. Let $M$ denote the maximal graded ideal of $R$.\sk

Let ${\frak m} = (\bar x_1,\ldots,\bar x_n)$ be the maximal
homogeneous ideal of $A$. We can refine the $\Bbb N$-graded
structure of $R$ by a bigrading with
\begin{eqnarray*} \bideg \bar x_i & = & (1,0),\ i = 1,\ldots,r,\\
\bideg f_jt  & = & (1,1),\ j = 1,\ldots,s. \end{eqnarray*} It is
easy to verify that if $z \in R$ is a bihomogeneous element with
$\bideg z = (a,b)$, then $\deg z = a$. There is also the natural
bigrading $\bideg \bar x_i = (1,0)$ and $\bideg f_jt = (0,1)$. But
we shall see that the first bigrading is more suitable for our
investigation. \sk

Let $E$ be any bigraded $R$-module. Then $E$ has the natural
$\Bbb Z$-graded structure $E_a = \bigoplus_{b \in \Bbb
Z}E_{(a,b)}$. In particular,  $H_M^i(E)$ is both a bigraded
$R$-module and a $\Bbb Z$-graded $R$-module with
$$H_M^i(E)_a = \bigoplus_{b \in \Bbb Z}H_M^i(E)_{(a,b)}.$$

Let $R_+$ denote the ideal of $R$ generated by the elements
$f_jt$. Let $G = \oplus_{n \ge 0}I^n/I^{n+1}$ be the associated
graded ring of $I$. To estimate $a_i(R)$ we consider the
following two short exact sequences of bigraded $R$-modules:
\begin{eqnarray*} 0 \To R_+ \To R \To A \To 0,\\
0 \To R_+(0,1) \To R \To G \To 0. \end{eqnarray*}

>From the above short exact sequences we obtain the following long
exact sequences of bigraded local cohomology modules:
\begin{eqnarray}
\cdots \To H_M^{i-1}(A)_{(a,b)} \To H_M^i(R_+)_{(a,b)}
\To H_M^i(R)_{(a,b)} \To H_M^i(A)_{(a,b)} \To \cdots  \\
\cdots \To H_M^{i-1}(G)_{(a,b)} \To H_M^i(R_+)_{(a,b+1)} \To
H_M^i(R)_{(a,b)} \To H_M^i(G)_{(a,b)} \To \cdots
\end{eqnarray}

These sequences allow us to study the vanishing of the bigraded
local cohomology modules of $R$ by means of those of $A$ and
$G$.\sk

\begin{Lemma} \label{vanish1}
For a fixed integer $a$ assume that there is an integer $b_0$ such
that \par \noindent {\rm (i)} $H_M^{i-1}(A)_{(a,b)} = 0$ for $b
\ge b_0$,\par \noindent {\rm (ii)} $H_M^i(G)_{(a,b)} = 0$ for $b
> b_0$. \par \noindent Then $H_M^i(R)_{(a,b)} = 0$ for $b \ge
b_0$.
\end{Lemma}

\begin{pf} The assumptions (i) and (ii) implies that for $b \ge b_0$,\par
\noindent (i') $H_M^i(R_+)_{(a,b)} \To H_M^i(R)_{(a,b)}$ is
injective,
\par \noindent (ii') $H_M^i(R_+)_{(a,b+1)} \To
H_M^i(R)_{(a,b)}$ is surjective.
\par \noindent Since $H_M^i(R)$ is an artinian $R$-module,
$H_M^i(R)_{(a,b)} = 0$ for $b$ large enough. Once we have
$H_M^i(R)_{(a,b+1)} = 0$ for some integer $b \ge b_0$, we can use
(i') and (ii') to deduce first that $H_M^i(R_+)_{(a,b+1)} = 0$
and then that $H_M^i(R)_{(a,b)} = 0$. This can be continued until
$b = b_0$. \end{pf}\sk

\begin{Proposition} \label{bound1}
$a_i(R) \le \max\{a_{i-1}(A),a_i(G)\}$.
\end{Proposition}

\begin{pf} Fix an arbitrary integer $a >
\max\{a_{i-1}(A),a_i(G)\}$. Since $a > a_{i-1}(A)$,
$H_M^{i-1}(A)_a = 0$. Therefore,
$$H_M^{i-1}(A)_{(a,b)} = 0\ \text{for all}\ b.$$
Since $a > a_i(G)$, $H_M^i(G)_a = 0$. Therefore,
$$H_M^i(G)_{(a,b)} = 0\ \text{for all}\ b.$$
By Lemma \ref{vanish1} we obtain $H_M^i(R)_{(a,b)} = 0$ for all
$b$. That implies $H_M^i(R)_a = 0$ or, equivalently, $a_i(R) \le
\max\{a_{i-1}(A),a_i(G)\}$.
\end{pf}\sk

On the other hand, there is the following relation between the
maximal shifts of the terms of the minimal free resolutions of
$A$ and $R$. \sk

\begin{Proposition} \label{shift} $b_i(R) \ge b_i(A)$.
\end{Proposition}

\begin{pf} We consider a minimal free resolution
${\Bbb F}: 0 \To F_l \To \cdots \To F_0 \To R$ of $R$ as a
bigraded module over the polynomial ring $S :=
k[X_1,\ldots,X_n,Y_1,\ldots,Y_s]$. Let ${\Bbb F}^*$ denote the
exact sequence:
$$0 \To \bigoplus_{a \in \Bbb Z}(F_l)_{(a,0)} \To
\cdots \To \bigoplus_{a \in \Bbb Z}(F_0)_{(a,0)} \To \bigoplus_{a
\in \Bbb Z}R_{(a,0)} = A$$ It is clear that ${\Bbb F}^*$ is a
free resolution of $A$ as a graded module over the polynomial
ring $B = k[X_1,\ldots,X_n]$. \par To estimate the shifts of the
twisted free modules of ${\Bbb F}^*$ we consider a twisted free
$S$-module $S(-c,-d)$. Then $(S(-c,-d))_{(a,0)} = S_{(a-c,-d)}$
is a direct sum of $s -d-1 \choose s-1$ copies of $B_{a-c+d} =
B(-(c-d))_a$, where we set ${s-d-1 \choose s-1} = 0$ if $d > 0$.
Therefore, $\oplus_a (S(-c,-d))_{(a,0)}$ is the direct sum of $s
-d-1 \choose s-1$ copies of $B(-(c-d))$. If $S(-c,-d)$ runs over
all twisted free modules of $F_i$, then $b_i(R) = \max\{c\} \ge
\max\{c-d\}.$ Since $\max \{c-d\}$ is the maximum shift of the
$i$th term of ${\Bbb F}^*$, we have $\max\{c-d\} \ge b_i(A)$. So
we obtain $b_i(R) \ge b_i(A)$. \end{pf} \sk

>From the above propositions we can easily derive upper and lower
bounds for $a^*(R)$ and $\reg(R)$ in terms of $A$ and $G$. \sk

\begin{Theorem} \label{a-invariant1} Let $s$ denote the minimal
number generators of $I$. Then $$a^*(A)-s \le a^*(R) \le
\max\{a^*(A),a^*(G)\}.$$
\end{Theorem}

\begin{pf} By definition we have $a^*(E) = \max\{a_i(E)|\ i \ge 0\}$
for any finitely generated graded $R$-module $E$. Therefore, from
Proposition \ref{bound1} we immediately obtain the upper bound
$a^*(R) \le \max\{a^*(A),a^*(G)\}.$ On the other hand, by
Corollary \ref{link2} we have
\begin{eqnarray*}
a^*(A) & = & \max\{b_i(A)|\ i\ge 0\}-n,\\
a^*(R) & = & \max\{b_i(R)|\ i \ge 0\} - n - s.
\end{eqnarray*}
Therefore, from Proposition \ref{shift} we can immediately deduce
the lower bound $a^*(A) - s \le a^*(R)$. \end{pf}

\begin{Theorem} \label{regularity1}
$\reg(A) \le \reg(R) \le \max\{\reg(A)+1,\reg(G)\}$
\end{Theorem}

\begin{pf} By Corollary \ref{link2} we have
$\reg(E) = \max\{a_i(E) + i|\ i \ge 0\}$ for any finitely
generated graded $R$-module $E$. By Proposition \ref{bound1},
$a_i(R)+i \le \max\{a_{i-1}(A)+i,a_i(G)+i\}$. Hence we get the
upper bound $\reg(R) \le \max\{\reg(A)+1,\reg(G)\}.$ On the other
hand, using Proposition \ref{shift} we obtain the lower bound
\begin{eqnarray*} \reg(A) &  = & \max\{b_i(A)+i|\ i \ge 0\}\\
& \le & \max\{b_i(R)+i|\ i \ge 0\}\ =\ \reg(R). \end{eqnarray*}
\end{pf} \sk

\begin{Corollary} Let $I$ be an ideal generated by a regular
sequence of $s$ forms of degree $c$. Then \par \noindent {\rm (i)
} $a^*(A)-s \le a^*(R) \le a^*(A) + s(c-1),$ \par \noindent {\rm
(ii) } $\reg(A) \le \reg(R) \le \max\{\reg(A)+1,\reg(A) +
s(c-1)\}$.
\end{Corollary}

\begin{pf} We have $G \cong (A/I)[z_1,\ldots,z_s]$, where
$z_1,\ldots,z_s$ are indeterminates. By Proposition
\ref{extension2}, this implies $a^*(G) = a^*(A/I) - s$ and
$\reg(G) = \reg(A/I)$. On the other hand, Proposition
\ref{quotient} gives $a^*(A/I) = a^*(A) + sc$ and $\reg(A/I) =
\reg(A) + s(c-1)$. Therefore, $a^*(G) = a^*(A)+s(c-1)$ and
$\reg(G) = \reg(A) + s(c-1)$. Hence the conclusion follows from
Theorem \ref{a-invariant1} and Theorem \ref{regularity1}.
\end{pf} \sk

The following example shows that the above upper and lower bounds
for $a^*(R)$ and $\reg(R)$ are sharp. \sk

\noindent{\bf Example.} Let $A = k[x_1,\ldots,x_n]$ and $I =
(x_1,\ldots,x_n)$. By Proposition \ref{extension2} we have
$a^*(A) = -n$ and $\reg(A) = 0$. In the next sections we shall
see that
\begin{eqnarray*} a^*(R) & = & \left\{ \begin{array}{lll} -2
& \text{for} & n = 1,\\ -n & \text{for} & n > 1. \end{array} \right.\\
\reg(R) & = & \left\{ \begin{array}{lll} 0 & \text{for} & n = 1,\\
1 & \text{for} & n > 1. \end{array} \right. \end{eqnarray*} \sk

\section{Local cohomology of Rees algebras of maximal graded
ideals}\sk

Let $A$ be a standard graded algebra over a field $k$. Let $\frak
m$ be the maximal graded ideal of $A$. From now on, $R$ will
denote the Rees algebra $A[{\frak m}t]$.\sk

We first note that $R$ has a bigraded automorphism $\psi$ induced
by the map $\psi(x) = xt$ and $\psi(xt) = x$ for any element $x
\in {\frak m}$. It is clear that if $f \in R$ is a bihomogeneous
element with $\bideg f = (a,b)$, then $\bideg \psi(f) = (a,a-b)$.
In particular, $\psi$ induces an isomorphism between $A$ and $G$
as $R$-modules. \par\sk

Since $A$ concentrates only in degree of the form $(a,0)$, we have
$$ H_M^i(A)_{(a,b)} =
\left\{ \begin{array}{lll}  0 & \text{for} & b \neq 0,\\
H^i_{\frak m}(A)_a  & \text{for} & b = 0. \end{array} \right.$$
>From this it follows that
$$ H_M^i(G)_{(a,b)} = \left\{ \begin{array}{lll} 0 & \text{for} & b \neq
a, \\ H^i_{\frak m}(A)_a & \text{for} & b = a. \end{array}
\right. $$ Therefore, using (1) and (2) we obtain
\begin{eqnarray} H_M^i(R_+)_{(a,b)} \cong H_M^i(R)_{(a,b)}\
\text{for}\ b \neq 0,\\
H_M^i(R_+)_{(a,b+1)} \cong H_M^i(R)_{(a,b)}\ \text{for}\ b \neq
a. \end{eqnarray}

The following lemma gives a complete description of the bigraded
local cohomology modules of $R$ in terms of those of $A$.\sk

\begin{Lemma} \label{vanish2} For any integer $a$ we have
$$H_M^i(R)_{(a,b)} = \left\{ \begin{array}{lll}
0 & \text{\rm if} & b \ge \max\{0,a+1\},\\
H_{\frak m}^i(A)_a & \text{\rm if} & 0 \le b < \max\{0,a+1\},\\
H_{\frak m}^{i-1}(A)_a & \text{\rm if} & \min\{0,a+1\} \le b < 0,\\
0 & \text{\rm if} & b < \min\{0,a+1\}. \end{array} \right.$$
\end{Lemma}

\begin{pf} By (3) and (4), $H_M^{i-1}(A)_{(a,b)} = 0$ for $b > 0$
and $H_M^i(G)_{(a,b)} = 0$ for $b \ge a+1$. Therefore, using
Lemma \ref{vanish1} we get $H_M^i(R)_{(a,b)} = 0$ for $b \ge
\max\{0,a+1\}$. By the above automorphism of $R$, this implies
$H_M^i(R)_{(a,b)} = 0$ for $a-b \ge \max\{0,a+1\}$ or,
equivalently, for $b < \min\{0,a+1\}$. \par

To prove that $H_M^i(R)_{(a,b)} = H_{\frak m}^i(A)_a$ for $0 \le
b < \max\{0,a+1\}$ we may assume that $a+1 > 0$. By (2) we have
the exact sequence
$$H_M^i(R_+)_{(a,a+1)} \To H_M^i(R)_{(a,a)}
\To H_M^i(G)_{(a,a)} \To H_M^{i+1}(R_+)_{(a,a+1)}.$$ By (5) we
have $H_M^j(R_+)_{(a,a+1)} \cong H_M^j(R)_{(a,a+1)} = 0$ for all
$j$. Therefore,
$$H_M^i(R)_{(a,a)} \cong H_M^i(G)_{(a,a)} = H_{\frak m}^i(A)_{a}.$$
Using (5) and (6) we obtain $H_M^i(R)_{(a,b)} \cong
H_M^i(R)_{(a,b+1)}$ for $0 \le b < a$. Thus, $H_M^i(R)_{(a,b)}
\cong H_{\frak m}^i(A)_a$ for $0 \le b < a+1 = \max\{0,a+1\}$.
\par

To prove that $H_M^i(R)_{(a,b)}  = H_{\frak m}^{i-1}(A)_a$ for
$\min\{0,a+1\} \le b < 0$ we may assume that $a+1 < 0$. By (1) we
have the exact sequence
$$H_M^{i-1}(R)_{(a,0)} \To H_M^{i-1}(A)_{(a,0)}
\To H_M^i(R_+)_{(a,0)} \To H_M^i(R)_{(a,0)}.$$

Since $a +1 < 0$, we have $H_M^j(R)_{(a,0)} = 0$ for all $j$.
Therefore,
$$H_M^i(R_+)_{(a,0)} \cong H_M^{i-1}(A)_{(a,0)} =
H_{\frak m}^{i-1}(A)_a.$$ By (2) we have the exact sequence
$$H_M^{i-1}(G)_{(a,-1)} \To H_M^i(R_+)_{(a,0)} \To
H_M^i(R)_{(a,-1)} \To H_M^i(G)_{(a,-1)}.$$ By (4) we have
$H_M^j(G)_{(a,-1)} = 0$ for all $j$. Therefore,
$$H_M^i(R)_{(a,-1)} \cong H_M^i(R_+)_{(a,0)} \cong
H_{\frak m}^{i-1}(A)_a.$$ Using (5) and (6) we obtain
$H_M^i(R)_{(a,b)} \cong H_M^i(R)_{(a,b+1)}$ for $a+1 \le b < -1$.
Thus, $H_M^i(R)_{(a,b)} \cong H_{\frak m}^{i-1}(A)_a$ for
$\min\{0,a+1\} = a+1 \le b \le -1$.
\end{pf} \sk

The above lemma on the bigraded local cohomology modules of $R$
can be formulated for the $\Bbb N$-graded structure as follows.
\sk

\begin{Theorem} \label{vanish3}
$$H_M^i(R)_a = \left\{ \begin{array}{lll}
\bigoplus_{a+1\, \text{\rm copies}}H_{\frak m}^i(A)_a & \text{\rm
if}
& a \ge 0,\\
0 & \text{\rm if} & a = -1,\\
\bigoplus_{-(a+1)\, \text{\rm copies}}H_{\frak m}^{i-1}(A)_a &
\text{\rm if} & a \le -2. \end{array} \right.$$
\end{Theorem}

\begin{pf} The statement follows from the formula
$$H_M^i(R)_a = \bigoplus_{b \in \Bbb Z}H_M^i(R)_{(a,b)}$$
and Lemma \ref{vanish2}. Indeed, if $a \ge 0$, $\max\{0,a+1\} =
a+1$ and $\min\{0,a+1\} = 0$. Therefore,
$$H_M^i(R)_{(a,b)} = \left\{ \begin{array}{lll}
0 & \text{\rm for} & b \ge a+1,\\
H_{\frak m}^i(A)_a & \text{\rm for} & 0 \le b \le a,\\
0 & \text{\rm for} & b < 0. \end{array} \right.$$ If $a = -1$,
$\max\{0,a+1\} =  \min\{0,a+1\} = 0$. Hence $H_M^i(R)_{(-1,b)} =
0$ for all $b$. If $a \le -2$, we have
$$H_M^i(R)_{(a,b)} = \left\{ \begin{array}{lll}
0 & \text{\rm for} & b \ge 0,\\
H_{\frak m}^{i-1}(A)_a & \text{\rm for} & a+1 \le b < 0,\\
0 & \text{\rm for} & b \le a. \end{array} \right.$$
 \end{pf} \sk

In the following we will denote $a_j(A)$ by $a_j$ for all $j$. An
immediate consequence of Theorem \ref{vanish3} is the following
formula for the depth of the Rees algebra (see [HM] for the depth
of the Rees algebra of an arbitrary ideal). \sk

\begin{Proposition} \label{depth} $$\depth R = \max\{i|\ H_{\frak
m}^j(A)_a = 0\ \text{\rm for}\ a \neq -1, j < i-1,\ \text{\rm
and}\ a_{i-1} < 0\}.$$ In particular, $\depth R \ge \depth A.$
\end{Proposition}

\begin{pf} We have $H_M^i(R) = 0$ if and only if
$H_{\frak m}^{i-1}(A)_a = 0$ for $a \le -2$ and $H_{\frak
m}^i(A)_a = 0$ for $a \ge 0$. Putting this in the formula
$$\depth R = \max\{i|\ H_M^j(R) = 0\ \text{\rm for}\ j < i\}.$$
we obtain the conclusion. \end{pf}\sk

If $\depth A = 0$, we must have $\depth R = 0$. If $\depth A = 1$,
we can not have $H_{\frak m}^1(A)_a = 0$ for $a \neq -1$. For, we
have $H_{\frak m}^1(A)_{-1} \cong H_\mm^0(A/xA)_0 = 0$, where $x$
is a regular linear form of $A$. Therefore, Proposition
\ref{depth} gives $\depth R = 1$ if $\depth A = 1$. If $\depth A
\ge 2$, $\depth R$ can be arbitrarily large than $\depth A$. \sk

\noindent{\bf Example.} Let $2 \le r < d$ be arbitrary positive
numbers. We will construct a graded algebra $A$ with $\depth A =
r$ and $\depth R = d+1$ (i.e. $R$ is a Cohen-Macaulay ring).\par
Let $T = k[x_1,\ldots,x_d]$ and $\nn$ the maximal graded ideal of
$T$. Let $E$ be the $r$th syzygy module of $k$ over $T$. Then
$H_\nn^i(E) = 0$ for $i \neq r,d$,
$$H_\nn^r(E)_a = \left\{ \begin{array}{lll}
0 & \text{\rm for} & a \ne 0,\\
k & \text{\rm for} & a = 0, \end{array} \right.$$ and $H_\nn^d(E)
= H_\nn^d(T)$. \par Let $C$ be the idealization of the graded
$T$-module $E(r-1)$ (see e.g. [N, p.2]). Since $E$ is generated by
elements of degree $r$, $E(r-1)$ is generated by elements of
degree 1. Hence $C$ is a standard graded algebra over $k$. By the
construction of the idealization we have a natural exact sequence
of the form $0 \to E(r-1) \to C \to T \To 0$, where all
homomorphisms have degree 0. Let $\mm_C$ denote the maximal
graded ideal of $C$. Then $H_{\mm_C}^i(C) = H_\nn^i(E(r-1))$ for
$i \neq d$ and there is the exact sequence $H_\nn^d(E(r-1)) \To
H_{\mm_C}^d(C) \To H_\nn^d(T) \To 0.$ From this it follows that
$H_{\mm_C}^i(C) = 0$ for $i \neq r,d$,
$$H_{\mm_C}^r(C)_a = \left\{ \begin{array}{lll}
0 & \text{\rm for} & a \ne 1-r,\\
k & \text{\rm for} & a = 1-r, \end{array} \right.$$ and $a_d(C) =
a_d(T) = -d$. \par Now let $A$ be the $(r-1)$th Veronese subring
of $C$. Then $H_\mm^i(A)_a = H_{\mm_C}^i(C)_{a(r-1)}$ [GW,
Theorem 3.1.1]. Hence $H_\mm^i(A) = 0$ for $i \neq r,d$,
$$H_\mm^r(A)_a = \left\{ \begin{array}{lll}
0 & \text{\rm for} & a \ne -1,\\
k & \text{\rm for} & a = -1, \end{array} \right.$$ and $a_d \le
-1$. Therefore, $\depth A = r$ and $\depth R = d+1$ by
Proposition \ref{depth}. \sk

\noindent
{\bf Remark.\ }The above example is inspired by Goto's
construction of Buchsbaum local rings of minimal multiplicity with
local cohomology modules of given lengths [G, Example (4.11)(2)].
Evans and Griffith [EvG] have constructed graded domains $A$ whose
local cohomology modules $H_\mm^i(A)$, $i <d$, are isomorphic to
given graded modules of finite length with a shifting. Since the
shift could not be computed explicitly, we can not use their
construction for our purpose.\sk

Despite the eventually big difference between the depths of a
given ring and its Rees algebra, the depth of the Rees algebra of
a polynomial extension $A[z]$ of $A$ is rather rigid. \sk

\begin{Corollary} \label{difference} Let $R^*$ denote the Rees algebra of a
polynomial ring $A[z]$ over $A$ in one variable. Put $s = \depth
A$. Then
$$\depth R^* = \left\{ \begin{array}{lll} s+1 & \text{\rm
if} & a_s \ge 0,\\ s+2 & \text{\rm if} & a_s < 0.
\end{array} \right.$$
\end{Corollary}

\begin{pf} By Proposition \ref{extension1}, $H_{\frak n}^i(A[z]) = 0$
if $H_{\frak m}^{i-1}(A) = 0$ and $H_{\frak n}^i(A[z])$ is a
module of infinite length with $a_i(A[z]) = a_{i-1}-1$ if
$H_{\frak m}^{i-1}(A) \neq 0$. As a consequence, $H_{\frak
n}^i(A[z])_a = 0\ \text{\rm for}\ a \neq -1$ if and only if
$H_{\frak m}^{i-1}(A) = 0$. Therefore,
\begin{eqnarray*} \depth R^* & = & \max\{i|\ H_{\frak n}^j(A[z])_a = 0\ \text{\rm
for}\ a \neq -1, j < i-1,\ \text{\rm and}\ a_{i-1} < 0\}\\ & = &
\max\{i|\ H_{\frak m}^j(A) = 0\ \text{\rm for}\ j < i-2,\
\text{\rm and}\ a_{i-2} < 0\}
\end{eqnarray*}
Since $s = \max\{i|\ H_{\frak m}^j(A) = 0\ \text{\rm for}\ j <
i\}$, we get
$$\max\{i|\ H_{\frak m}^j(A) = 0\ \text{\rm for}\ j < i-2,\
\text{\rm and}\ a_{i-2} < 0\} =
\left\{ \begin{array}{lll} s+1 & \text{\rm if} & a_s \ge 0,\\
s+2 & \text{\rm if} & a_s < 0.
\end{array} \right. $$ \end{pf} \sk

The following criterion for the Cohen-Macaulayness of $R$ can be
also derived from a more general criterion for the
Cohen-Macaulayness of the Rees algebra of an arbitrary ideal of
Trung and Ikeda [TI]. \sk

\begin{Corollary} \label{CM} $R$ is a Cohen-Macaulay ring if and
only if $H_{\frak m}^j(A)_a = 0$ for $a \neq -1$, $j < d$, and
$a_d < 0$. \end{Corollary}

\begin{pf} This follows from Proposition \ref{depth} (the case
$\depth R = d+1$). \end{pf}\sk

Another consequence of Theorem \ref{vanish3} is the following
formula for $a_i(R)$. This formula is crucial for the estimation
of $a^*(R)$ and $\reg(R)$ in the next section. \sk

\begin{Proposition} \label{bound2}
$$a_i(R) = \left\{ \begin{array}{lll} a_i & \text{\rm if} & a_i \ge
0,\\ \max\{a|\ a \le -2\ \text{\rm and}\ H_{\frak m}^{i-1}(A)_a
\neq 0\} & \text{\rm if} & a_i < 0. \end{array} \right.$$ In
particular, $a_i(R) = a_{i-1}$ if $a_{i-1} \le -2$ and $a_i < 0$.
\end{Proposition}

\begin{pf} Note that $H_{\frak m}^j(A)_a = 0$ for $a > a_j$
and $H_{\frak m}^i(A)_{a_j} \neq 0$ if $a_j \neq -\infty$. If
$a_i \ge 0$, using Theorem \ref{vanish3} we get $H_M^i(R)_a = 0$
for $a > a_i$ and $H_M^i(R)_{a_i} \neq 0$, hence $a_i(R) = a_i$.
If $a_i < 0$, we get $H_M^i(A)_a = 0$ for $a \ge -1$. For $a \le
-2$, $H_M^i(R)_a \neq 0$ if and only if $H_{\frak m}^{i-1}(A)_a
\neq 0$. Therefore, $a_i(R) = \max\{a|\ a \le -2\ \text{\rm and}\
H_{\frak m}^{i-1}(A)_a \neq 0\}$, which is exactly $a_{i-1}$ if
$a_{i-1} \le -2$. \end{pf} \sk

From Proposition \ref{bound2} we immediately obtain the following
bounds for $a_i(R)$.\sk

\begin{Corollary}
\noindent {\rm (i)} $a_i(R) \le \max\{a_{i-1},a_i\}$ for $i \le
d$,\par \noindent {\rm (ii)} $a_{d+1}(R) \le \min\{-2,a_d\}$.
\end{Corollary} \sk

\noindent{\bf Example. } Let $A = k[x_1,\ldots,x_n]$. We know
that $a_i = -\infty$ for $i \neq n$ and $a_n = -n$ with
$H_\mm^n(A)_a \neq 0$ for $a \le -n$. Therefore, $a_i(R) =
-\infty$ for $i \neq n+1$ and
$$a_{n+1}(R) = \left\{ \begin{array}{lll} -2 & \text{\rm if} &
n = 1,\\ -n & \text{\rm if} & n > 1. \end{array} \right.$$ \sk

One may expect that $a_i(R) = -2$ if $a_{i-1} > -2$ and $a_i <
0$. But the following example shows that is not always the case.
\sk

\noindent{\bf Example.} Let $\Delta$ be the simplicial complex on
ten vertices $\{1,\ldots,10\}$ with the maximal faces
$$\{1,2,6\},\{2,6,7\},\{2,3,7\},\{3,7,8\},\{3,4,8\},\{4,5,8\},\{1,4,5\},\{9,10\}.$$
Note that $\Delta$ is topologically the disjoint union of a
circle and a point.

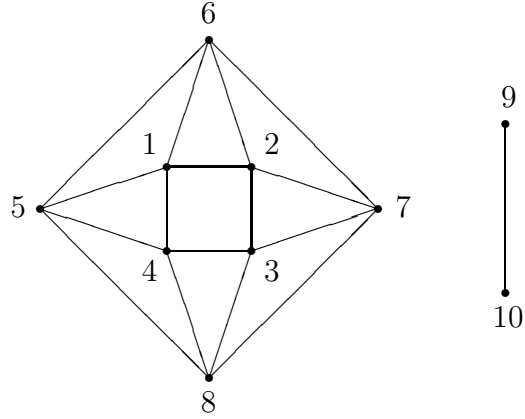
\begin{picture}(300,170)(-100,0) \setlength{\unitlength}{16pt}
\put(1,5){\circle*{.2}}\put(4,4){\circle*{.2}}
\put(4,6){\circle*{.2}} \put(5,1){\circle*{.2}}
\put(5,9){\circle*{.2}} \put(6,4){\circle*{.2}}
\put(6,6){\circle*{.2}} \put(9,5){\circle*{.2}}
\put(12,3){\circle*{.2}}\put(12,7){\circle*{.2}}
\put(4,6){\line(0,-1){2}} \put(4,6){\line(1,0){2}}
\put(6,6){\line(0,-1){2}} \put(4,4){\line(1,0){2}}
\put(1,5){\line(3,-1){3}} \put(1,5){\line(3,1){3}}
\put(4,4){\line(1,-3){1}} \put(4,6){\line(1,3){1}}
\put(5,1){\line(1,3){1}} \put(5,9){\line(1,-3){1}}
\put(6,4){\line(3,1){3}} \put(6,6){\line(3,-1){3}}
\put(1,5){\line(1,1){4}} \put(1,5){\line(1,-1){4}}
\put(5,1){\line(1,1){4}} \put(5,9){\line(1,-1){4}}
\put(12,3){\line(0,1){4}} \put(3.4,6.3){1} \put(6.3,6.3){2}
\put(6.3,3.3){3} \put(3.4,3.3){4} \put(0.3,4.8){5}
\put(4.8,9.4){6} \put(9.4,4.8){7} \put(4.8,0.2){8}
\put(11.7,2.2){10} \put(11.9,7.4){9}
\end{picture}

Let $T = k[x_1,\ldots,x_{10}]$ and $\nn$ the maximal graded ideal
of $T$. Let $I$ be the monomial ideal of $\Delta$ in $T$:
\begin{eqnarray*} I & = & (x_2x_{10}, x_2x_9, x_2x_8, x_2x_5,
x_2x_4, x_3x_{10},x_3x_9, x_3x_5, \\ & & x_1x_3, x_5x_{10},
x_4x_{10}, x_1x_{10}, x_5x_9, x_4x_9, x_1x_9, x_1x_8, x_8x_9,
x_8x_{10}, \\ & & x_1x_7, x_4x_7, x_5x_7, x_7x_9, x_7x_{10},
x_3x_6, x_4x_6, x_5x_6, x_6x_8, x_6x_9, x_6x_{10}).
\end{eqnarray*} By Hochster's formula for the local cohomology
modules of the Stanley-Reisner ring $k[\Delta]= T/I$ (see e.g.
[BH, Theorem 5.8]) we have $H_\nn^0(k[\Delta]) = 0$,
$$H_\nn^1(k[\Delta])_a = \left\{ \begin{array}{lll}
0 & \text{\rm for} & a \neq 0,\\
k & \text{\rm for} & a = 0, \end{array} \right.$$
$$H_\nn^2(k[\Delta])_a = \left\{ \begin{array}{lll}
0 & \text{\rm for} & a > 0\ \text{and}\ a=-1,\\
k & \text{\rm for} & a = 0\ \text{and}\ a \le -2, \end{array}
\right.$$ and $H_\nn^3(k[\Delta])_a = 0$  for $a > -2$, while
$H_\nn^3(k[\Delta])_a \neq 0$ for $a \le -2$. From the short exact
sequence $0 \To I \To T \To k[\Delta] \To 0$ we get $H_\nn^i(I) =
0$ for $i \neq 2,3,4,10$,
$$H_\nn^2(I)_a = \left\{ \begin{array}{lll}
0 & \text{\rm for} & a \neq 0,\\
k & \text{\rm for} & a = 0, \end{array} \right.$$
$$H_\nn^3(I)_a = \left\{ \begin{array}{lll}
0 & \text{\rm for} & a > 0\ \text{and}\ a=-1,\\
k & \text{\rm for} & a = 0\ \text{and}\ a \le -2, \end{array}
\right.$$ and $H_\nn^4(I)_a = 0$  for $a > -2$, $H_\nn^4(I)_a
\neq 0$ for $a \le -2$, $H_\nn^{10}(I)_a = 0$ for $a > -10$,
$H_\nn^{10}(I)_a \neq 0$ for $a \le -10$. \par Let $A$ be the
idealization of the graded $T$-module $I(1)$ (see e.g. [N, p.2]).
Since $I(1)$ is generated by elements of degree 1, $A$ is a
standard graded algebra over $k$. By the construction of the
idealization we have a natural exact sequence of the form $0 \To
I(1) \To A \To T \To 0$, where all homomorphisms are of degree 0.
This exact sequence yields $H_\mm^i(A) = H_\nn^i(I(1))$ for $i
\neq 10$ and the exact sequence $H_\nn^d(I(1)) \To H_\mm^d(A) \To
H_\nn^d(T) \To 0.$ Using the above formula for $H_\nn^i(I)$ we
get $H_\mm^i(A) = 0$ for $i\neq 2,3,4,10$,
$$H_\mm^2(A)_a = \left\{ \begin{array}{lll}
0 & \text{\rm for} & a \neq -1,\\
k & \text{\rm for} & a = -1, \end{array} \right.$$
$$H_\nn^3(A)_a = \left\{ \begin{array}{lll}
0 & \text{\rm for} & a > -1 \ \text{and}\ a=-2,\\
k & \text{\rm for} & a = -1 \ \text{and}\ a \le -3, \end{array}
\right.$$ and $H_\mm^4(A)_a = 0$  for $a > -3$, $H_\mm^4(A)_a \neq
0$ for $a \le -3$, $H_\mm^{10}(A)_a = 0$ for $a
> -10$, $H_\nn^{10}(A)_{-10} \neq 0$. In particular,
$a_i = -\infty$ for $i \neq 2,3,4,10$, $a_2  = a_3 = -1$, $a_4 =
-3$, $a_{10} = -10$. Applying Proposition \ref{bound2} we obtain
$a_i(R) = -\infty$ for $i \neq 4,5,11$ and $a_4(R) = a_5(R) = -3$,
$a_{11}(R) = -10$. In particular, $a_4(R) = -3$ though $a_3 = -1$
and $a_4 = -3 < 0$. \sk

\section{Regularity of Rees algebras of maximal graded ideals}\smallskip

As in the last section, let $R = A[{\frak m}t]$ be the Rees
algebra of the maximal graded ideal $\frak m$ of a standard
graded algebra $A$ over a field. The goal of this section is to
estimate $a^*(R)$ and $\reg(R)$ in terms of $a^*(A)$ and
$\reg(A)$. \sk

Let $n$ be the embedding dimension of $A$, that is, $n =
\dim_kA_1$. We can consider $A$ as a module over the polynomial
ring $B = k[x_1,\ldots,x_n]$ and $R$ as a module over the
polynomial ring $S = k[x_1,\ldots,x_n,y_1,\ldots,y_n]$.\sk

We have the following relationships between the invariants
$a_j^*(A)$ and $a_j^*(R)$. \sk

\begin{Proposition} \label{a-invariant2}
For any integer $j \ge 0$ we have \par \noindent {\rm (i)}
$a_{j-n}^*(A)-n \le a_j^*(R) \le a_j^*(A)$. \par \noindent{\rm
(ii)} $a_j^*(R) = a_j^*(A)$ if and only if $a_j^*(A) \ge 0$ or
$a_{j-1}^*(A) = a_j^*(A) \le -2$.
\end{Proposition}

\begin{pf} By Theorem \ref{link1}(i) and Proposition \ref{shift} we have
$$a_j^*(R) = b_{2n-j}^*(R) - 2n \ge b_{2n-j}^*(A) - 2n =
a_{j-n}^*(A)-n.$$ Since $G \cong A$, Theorem \ref{regularity1}
implies $a_i(R) \le \max\{a_{i-1}(A),a_i(A)\}$ for all $i$. From
this it follows that $a_j^*(R) \le a_j^*(A)$. So we obtain (i).
\par To prove (ii) we choose $i \le j$ such that $a^*_j(A) =
a_i$. \par If $a_j^*(A) \ge 0$, then $a_i(R) = a_i$ by Corollary
\ref{bound2}. Hence $a^*_j(R) \ge a_i(R) = a_j^*(A)$. By (i) this
implies $a_j^*(R) = a_j^*(A)$.
\par If $a_j^*(A) = -1$, $a_i(R) \le -2$ for all $i \le j$ by
Proposition \ref{bound2}. Hence $a_j^*(R) \le -2 < a_j^*(A)$. \par
If $a_{j-1}^*(A) = a_j^*(A) \le -2$, we may assume that $i < j$.
Then $a_{i+1} \le a_i \le -2$. Hence $a_{i+1}(R) = a_i$ by
Proposition \ref{bound2}. Thus, $a_j^*(R) \ge a_{i+1}(R) =
a_j^*(A)$. By (i) this implies $a_j^*(R) = a_j^*(A)$. \par If
$a_{j-1}^*(A) < a_j^*(A) \le -2$, we have $a_{j-1} < a_j =
a_j^*(A) \le -2$. By Proposition \ref{bound2}, this implies
$a_j(R) = a_{j-1} < a_j^*(A)$. By (i), $a_{j-1}^*(R) \le
a_{j-1}^*(A)$. Therefore, $a_j^*(R) = \max\{a_{j-1}^*(R),a_j(R)\}
< a_j^*(A).$ So we have proved (ii).
\end{pf} \sk

Proposition \ref{a-invariant2} can be formulated in terms of the
maximal shifts of the minimal free resolution of $A$ as
follows.\sk

\begin{Corollary} For any integer $j \ge 0$ we have \par
\noindent {\rm (i)} $b_j^*(A) \le b_j^*(R) \le
b_{j-n}^*(A)+n$,\par \noindent {\rm (ii)} $b_j^*(R) =
b_{j-n}^*(A)+n$ if and only if $b_{j-n}^*(A) \ge n$ or
$b_{j-n+1}^*(A) = b_{j-n}^*(A)\le n-2.$
\end{Corollary}

\begin{pf} By Theorem \ref{link1}(i) we have
$$b_j^*(A) = a_{n-j}^*(A)+n,\ b_j^*(R) = a_{2n-j}^*(R)+2n,\
b_{j-n}^*(A) = a_{2n-j}^*(A)+n.$$ Therefore, the conclusion
follows from Proposition \ref{a-invariant2}. \end{pf}\sk

\begin{Theorem} \label{a-invariant3} {\rm (i)} $a^*(A)-n \le
a^*(R)\le a^*(A)$. \par \noindent{\rm (ii)} $a^*(R) = a^*(A)$ if
and only if $a^*(A) \neq -1$.
\end{Theorem}

\begin{pf} From Proposition \ref{a-invariant2}(i) we immediately obtain (i).
To prove (ii) we assume first that $a^*(A) \neq -1$. Choose $j$
such that $a^*(A) = a_{j-1}^*(A) = a_j^*(A)$. Then $a_j^*(R) =
a_j^*(A)$ by Proposition \ref{a-invariant2}(ii). It follows that
$$a^*(R) \ge a_j^*(A) = a^*(A).$$
By (i) this implies $a^*(R) = a^*(A)$. Now assume that $a^*(A) =
-1$. Then $a_j^*(A) \le -1$ for all $j$. If $a_j^*(A) \le -2$, we
have $a_j^*(R) \le a_j^*(A) \le -2$ by Theorem
\ref{a-invariant2}(i). If $a_j^*(A) = -1$, then Theorem
\ref{a-invariant2}(ii) implies $a_j^*(R) < a_j^*(A)$. Thus,
$a^*(R) \le -2$. So we have proved that $a^*(R) = a^*(A)$ if and
only if $a^*(A) \neq -1$. \end{pf} \sk

\noindent{\bf Example.} Let $A = k[x_1,\ldots,x_n]$. From the
formula for $a_i(A)$ and $a_i(R)$ in the last section we get
$a^*(A) = -n$ while
$$a^*(R) =  \left\{ \begin{array}{lll} -2
& \text{for} & n = 1,\\ -n & \text{for} & n > 1. \end{array}
\right.$$ \sk

If $a^*(A) = -1$, we may expect that $a^*(R) = -2$. But that is
not always the case. For instance, in the last example of Section
4 we have $a^*(A) = -1$ and $a^*(R) = -3$. \sk

We may also formulate Theorem \ref{a-invariant3} in terms of the
maximal shifts of the minimal free resolution of $A$ as follows.
\sk

\begin{Corollary} {\rm (i)} $b^*(A) \le b^*(R) \le b^*(A)+n.$ \par
\noindent {\rm (ii)} $b^*(R) = b^*(A)+n$ if and only if $b^*(A)
\ne n-1.$
\end{Corollary}

\begin{pf} By Corollary \ref{link2}(i) we have $b^*(A) = a^*(A)+n$
and $b^*(R) = a^*(R)+2n$. Hence the conclusion follows from
Theorem \ref{a-invariant3}. \end{pf} \sk

Now we study the relationships between the partial regularities of
a given graded algebra $A$ and its Rees algebra $R$. \sk

\begin{Proposition} \label{regularity2} For any integer $j \ge 0$
we have \par \noindent {\rm (i)} $\reg_{j-n}(A) \le \reg_j(R) \le
\reg_j(A)+1$. \par \noindent {\rm (ii)} $\reg_j(R) = \reg_j(A)+1$
if and only if there is an integer $i <j$ such that $\reg_j(A) =
a_i+i$ and $a_i\le -2$. \end{Proposition}

\begin{pf} By Theorem \ref{link1}(ii) and Proposition \ref{shift} we have
$$\reg_j(R) = (2n-j)\text{-}\reg(R) \ge (2n-j)\text{-}\reg(A) =
\reg_{j-n}(A).$$ Since $G \cong A$, Theorem \ref{regularity1}
implies $a_i(R) \le \max\{a_{i-1}(A),a_i(A)\}$ for all $i$. From
this it follows that $\reg_j(R) \le \reg_j(A)+1$. So we obtain
(i). \par To prove (ii) we assume first that $\reg_j(R)=
\reg_j(A)+1$. Let $i \le j$ be an integer such that $\reg_j(R) =
a_i(R)+i$. Then $a_i(R) = \reg_j(A)-i+1$. On the other hand,
$a_i(R)\le \max\{a_{i-1},a_i\}$ by Proposition \ref{bound1}.
Since $a_{i-1} \le \reg_j(A)-i+1$ and $a_i \le \reg_j(A)-i$, we
must have $a_i(R) = a_{i-1}$. This implies $\reg_j(A) =
a_{i-1}+i-1$ and, by Proposition \ref{bound2}, $a_{i-1} \le -2$.
\par Conversely, assume that there is an integer $i < j$ such
that $\reg_j(A) = a_i+i$ and $a_i\le -2$. Then $a_{i+1} + i + 1
\le a_i + i$. Therefore, $a_{i+1} \le a_i-1 < 0$. By Corollary
\ref{bound2}, we get $a_{i+1}(R) = a_i$. Hence
$$\reg_j(R) \ge a_{i+1}(R) + i +1 = a_i + i + 1 = \reg_j(A)+1.$$
By (i) this implies $\reg_j(R) = \reg_j(A)$. \end{pf} \sk

We may formulate Proposition \ref{regularity2}(i) for the partial
regularity $j\text{-}\reg(R)$ of Bayer, Charalambous, and Popescu
(see Section 1). But, unlike the estimation for $b_j^*(R)$, we
are not able to express the condition of Proposition
\ref{regularity2}(ii) in terms of the maximal shifts of the
minimal free resolution of $A$. \sk

\begin{Corollary} For any integer $j \ge 0$ we have \par
\noindent {\rm (i)} $j\text{-}\reg(A) \le j\text{-}\reg(R) \le
(j-n)\text{-}\reg(A)+1$. \par \noindent {\rm (ii)}
$j\text{-}\reg(R) = (j-n)\text{-}\reg(A)+1$ if
$(j-n)\text{-}\reg(A) = b_i(A)-i$ and $b_i(A)\le n-2$ for some
index $i > j-n$ at which $A$ has an extremal Betti number.
\end{Corollary}

\begin{pf} By Theorem \ref{link1}(ii) we have
$$j\text{-}\reg(A) = \reg_{n-j}(A),\ j\text{-}\reg(R) = \reg_{2n-j}(R),\
(j-n)\text{-}\reg(A) = \reg_{2n-j}(A).$$ Therefore, (i) follows
from Proposition \ref{regularity2}(i). For (ii) we have $b_i(A) =
a_{n-i}+n$ by Corollary \ref{link3}, hence $\reg_{2n-j}(A) =
a_{n-i}+n-i$ and $a_{n-i} \le -2$. By Theorem
\ref{regularity2}(ii) this implies $\reg_{2n-j}(R) =
\reg_{2n-j}(A)+1$. Thus, $j\text{-}\reg(R) =
(j-n)\text{-}\reg(A)+1$. \end{pf} \sk

The following result which is an immediate consequence of Proposition \ref{regularity2} 
gives precise information on the value of the
regularity of the Rees algebra of the maximal graded ideal. \sk

\begin{Theorem} \label{regularity3}
{\rm (i)} $\reg(A) \le \reg(R) \le \reg(A)+1$. \par \noindent
{\rm (ii)} $\reg(R) = \reg(A)+1$ if and only if there is an
integer $i$ such that $\reg(A) = a_i+i$ and $a_i \le -2$.
\end{Theorem}

\noindent{\bf Example.} Let $A = k[x_1,\ldots,x_n]$. We have
$\reg(A) = 0$ with $a_i = -\infty$ for $i < n$ and $a_n = -n$.
Therefore, $$\reg(R) = \left\{ \begin{array}{lll} 0 & \text{for} & n = 1,\\
1 & \text{for} & n > 1. \end{array} \right.$$ \sk

Let $T = A[z_1,\ldots,z_s]$ be a polynomial ring over $A$. Let
$R_s$ denote the Rees algebra of $T$ with respect to the maximal
graded ideal. It is well known that $\reg(T) = \reg(A)$. However,
the regularities of the Rees algebras $R_s$ and $R$ need not to
be the same. In fact, it may happen that $\reg(R) = \reg(A)$ but
$\reg(R_s) = \reg(A)+1$. \sk

\begin{Corollary} Let $i = \max\{j|\ \reg(A) = a_j+j\}$.
Put $c = \max\{0,a_i -2\}$. Then
$$\reg(R_s) = \left\{ \begin{array}{lll} \reg(A) & \text{for} & s < c,\\
\reg(A)+1 & \text{for} & s \ge c. \end{array} \right.$$
\end{Corollary}

\begin{pf} By Proposition \ref{extension2}
we have $a_{j+s}(T) = a_j - s$ for all $j \ge 0$. Since $\reg(T)
= \reg(A)$, we have $\reg(T) = a_{j+s}(T)+ (j+s)$ if and only if
$\reg(A) = a_j +j$. If $s < c$, then $a_j - s
> a_i - c \ge -2$ for all $j$ with $\reg(A) = a_j+j$.
>From this it follows that $a_{j+s}(T) > -2$ for all $j$ with
$\reg(T) = a_{j+s}(T)+(j+s)$. If $s \ge c$, then $\reg(T) =
a_{i+s}(T)+(i+s)$ with $a_{i+s}(T) = a_i -s \le -2$. Now we only
need to apply Theorem \ref{regularity3}(ii) to get the conclusion.
\end{pf}\sk

There is the following sufficient condition for the equality
$\reg(R) = \reg(A)+1$ in terms of the maximal shifts of the
minimal free resolution of $A$. \sk

\begin{Corollary} \label{sufficient} $\reg(R) = \reg(A)+1$ if
$\reg(A) = b_i(A)-i$ and $b_i(A) \le n-2$ for some index $i$ at
which $A$ has an extremal Betti number. \end{Corollary}

\begin{pf} By Corollary \ref{link3} we have $b_i(A) =
a_{n-i}+n$, hence $\reg(A) = a_{n-i}(A) + n-i$ and $a_{n-i} \le
-2$. By Proposition \ref{regularity2} this implies $\reg(R) =
\reg(A)+1$. \end{pf} \sk

Corollary \ref{sufficient} is not a necessary condition for the
equality $\reg(R) = \reg(A)+1$. \sk

\noindent{\bf Example.} Let $A=k[\Delta]$, where $\Delta$ is the simplicial complex on 7 
vertices $\{1,2,3,4,\\ 
5,6,7\}$ with the maximal faces $\{1,2,3\}$, $\{4,5,6\}$ and $\{5,6,7\}$. 
Using Hochster's formula for the local cohomology modules of $k[\Delta]$ (see e.g.\ [BH, 
Theorem 5.8]) we get $H^i_\mm(A)=0$ for $i\neq 1,3$, $H^1_\mm(A)_a=0$ for $a\neq 0$, 
$H^3_\mm(A)_a\neq 0$ for $a\leq -2$ and $H^3_\mm(A)_a=0$ for $a>-2$. From this it follows 
that $a_i(A)=-\infty$ for $i\neq 1,3$, $a_1(A)=0$ and $a_3(A)=-2$. Therefore, $\reg(A)=1$ 
and, by Theorem \ref{regularity3}, $\reg(R)=2$. On the other hand, $A$ has a $2$-linear 
$R$-resolution since $\reg(A)=1$. Hence the condition that $\reg(A)=b_i(A)-i$, and that 
$A$ has an extremal Betti number at $i$ is satisfied for $i=6$. But we have 
$b_6(A)=7>n-2=5$.
\sk

As applications, we will study the regularity of the Rees algebra
$R_{\ini} := B[{\mm}_{\ini}t]$ of the algebra $A_{\ini} :=
B/(\ini(J))$, where $\ini(J)$ denotes the initial ideal of $J$
with respect to an arbitrary term order and $\mm_{\ini}$ is the
maximal graded ideal of $A_{\ini}$.  \sk

\begin{Proposition} \label{ini3} {\rm (i)} $a_i(R) \le
a_i(R_{\ini})$ for all $i \ge 0$,\par \noindent {\rm (ii)} $a^*(R)
\le a^*(R_{\ini})$,\par \noindent {\rm (iii)} $\reg(R) \le
\reg(R_{\ini})$.
\end{Proposition}

\begin{pf} We only need to prove (i). Note that $a_i \le a_i(A_{\ini})$ by
[Sb, Theorem 3.3]. If $a_i \ge 0$, then $a_i(A_{\ini}) \ge 0$.
Applying Proposition \ref{bound2} we obtain
$$a_i(R) = a_i \le a_i(A_{\ini} ) = a_i(R_{\ini}).$$
If $a_i < 0$ and $a_i(A_{\ini} ) < 0$, then
\begin{eqnarray*} a_i(R) & = &
\max\{a|\ a\le -2\ \text{\rm and}\ H_{\frak m}^{i-1}(A)_a \neq 0\}\\
& \le & \max\{a|\ a\le -2\ \text{\rm and}\ H_{\frak
m}^{i-1}(A_{\ini})_a \neq 0\}\ = a_i(R_{\ini}). \end{eqnarray*} If
$a_i < 0$ and $a_i(A_{\ini} ) \ge 0$, we apply Proposition
\ref{bound2} again to see that
$$a_i(R) \le -2 < a_i(A_{\ini} ) = a_i(R_{\ini}).$$
 \end{pf}\sk

For the generic initial ideal $\Gin(J)$ of $I$ with respect to
the reverse lexicographic term order, we set $A_{\Gin} :=
B/\Gin(J)$ and $R_{\Gin} := B[\mm_{\Gin}t]$, where $\mm_{\Gin}$
is the maximal graded ideal of $A_{\Gin}$. The following results
show that in this case, $a^*(R_{\Gin})$ and $\reg(R_{\Gin})$
share the same lower and upper bounds of $a^*(R)$ and $\reg(R)$
as in Theorem \ref{a-invariant3} and Theorem \ref{regularity3}. In
particular, $\reg(R_{\Gin})$ differs from $\reg(R)$ at most by 1.
\sk

\begin{Proposition} {\rm (i)} $a^*(A)-n \le a^*(R_{\Gin}) \le
a^*(A)$,\par \noindent {\rm (ii)} $a^*(R_{\Gin}) = a^*(A)$ if and
only if $a^*(A) \neq -1$. \end{Proposition}

\begin{pf} By [T3, Corollary 1.4] we have $a^*(B/\Gin(J)) = a^*(A)$.
Therefore, the conclusion follows from Theorem \ref{a-invariant3}.
\end{pf}\sk

\begin{Proposition} {\rm (i)} $\reg(A) \le \reg(R_{\Gin}) \le
\reg(A)+1$, \par \noindent {\rm (ii)} $\reg(R_{\Gin}) = \reg(A)+1$
if there is an integer $i$ such that $\reg(A) = a_i+i$ and
$a_i(A) \le -2$.
\end{Proposition}

\begin{pf} By [T3, Corollary 1.4] we have
$\reg(A_{\Gin}) = \reg(A)$. Therefore, (i) follows from Theorem
4.3(i). By Theorem \ref{regularity3}(ii), the condition of (ii)
implies that $\reg(R) = \reg(A)+1$. By Proposition \ref{ini3}, we
have $\reg(R) \le \reg(R_{\Gin}) \le \reg(A)+1$. Therefore,
$\reg(R_{\Gin}) = \reg(A)+1.$ \end{pf} \sk

\noindent{\it Remark.} In spite of the above results one may ask
whether $a^*(R) = a^*(R_{\Gin})$ and $\reg(R) = \reg(R_{\Gin})$
always hold. Unfortunately, we were unable to settle this
question. \sk

\section*{References}\sk

\noindent [BM] D. Bayer and D. Mumford, What can be computed in
algebraic geometry? In "Computational Algebraic Geometry and
Commutative Algebra, Proceedings, Cortona, 1991" (D. Eisenbud and
L. Robbiano, Eds.), pp. 1-48, Cambridge Univ. Press, 1993. \par

\noindent [BCP] D.~Bayer, H.~Charalambous and S. Popescu,
Extremal Betti numbers and applications to monomial ideals, J.
Algebra 221 (1999), 497-512.\par

\noindent [BH] W. Bruns and J. Herzog, Cohen-Macaulay rings,
Cambridge, 1998.\par

\noindent [EiG] D.~Eisenbud and S.~Goto, Linear free resolutions
and minimal multiplicities, J.~Algebra~88 (1984), 89-133.\par

\noindent [EvG] E.~G. Evans and P.~.A. Griffith, Local cohomology
modules for normal domains, J. London Math. Soc. 19 (1979),
277-284. \par

\noindent [G] S. Goto, On the associated graded rings of
parameter ideals in Buchsbaum rings, J. Algebra 85 (1980),
490-534. \par

\noindent [GW] S. Goto and K. Watanabe, On graded rings I, J.
Math. Soc. Japan 30 (1978), 179-213. \par

\noindent [HM] S. Huckaba and T. Marley, Depth formulas for
certain graded rings associated to an ideal, Nagoya Math. J. 133
(1994), 57-69. \par

\noindent [JK] B. Johnston and D. Katz, Castelnuovo regularity
and graded rings associated to an ideal. Proc. Amer. Math. Soc.
123 (1995), no. 3, 727--734.\par

\noindent [N] M.~Nagata, Local rings, New York-London, 1962.
\par

\noindent [Sb] E. Sbarra, Upper bounds for local cohomology
modules of rings with a given Hilbert function, Thesis,
University of Essen, 2000.\par

\noindent [Sh] R.~Y. Sharp, Bass numbers in the graded case,
$a$-invariant formulas, and an analogue of Faltings' annihilator
theorem. J. Algebra 222 (1999), 246-270. \par

\noindent [T1] N.~V.~Trung, The Castelnuovo regularity of the
Rees algebra and the associated graded ring, Trans. Amer. Math.
Soc. 350 (1998), 2813-2832. \par

\noindent [T2] N.~V.~Trung, The largest non-vanishing degree of
graded local cohomology modules, J. Algebra 215 (1999), 481-499.
\par

\noindent [T3] N.~V.~Trung, Gr\"obner bases, local cohomology and
reduction number, to appear in Proc. Amer. Math. Soc.\par

\noindent [TI] N.~V.~Trung and S. Ikeda, When is the Rees algebra
Cohen-Macaulay?,  Comm. Algebra~17 (12) (1989), 2893-2922.

\end{document}